\documentclass[12pt,a4paper]{amsart}
\usepackage{amsmath}
\usepackage{amscd}

\newcommand{\N}{\mathbb{N}}
\newcommand{\Z}{\mathbb{Z}}
\newcommand{\cA}{\mathcal{A}}
\newcommand{\cC}{\mathcal{C}}
\newcommand{\cD}{\mathcal{D}}
\newcommand{\cF}{\mathcal{F}}
\newcommand{\cH}{\mathcal{H}}
\newcommand{\cI}{\mathcal{I}}
\newcommand{\cT}{\mathcal{T}}
\DeclareMathOperator{\End}{End}
\DeclareMathOperator{\Hom}{Hom}
\DeclareMathOperator{\Rep}{Rep}
\DeclareMathOperator{\Ob}{Ob}
\newtheorem{example}{Example}[section]
\newtheorem{lemma}[example]{Lemma}
\newtheorem{definition}[example]{Definition}
\newtheorem{proposition}[example]{Proposition}

\begin{document}
\title{Representation theory for dilute lattice models}
\author{Bruce W. Westbury}
\address{Mathematics Institute\\
University of Warwick\\
Coventry CV4 7AL}
\email{bww@maths.warwick.ac.uk}
\date{\today}
\maketitle
\tableofcontents

\section{Introduction}
In this paper we study the representation theory associated to
dilute lattice models. Given a solvable lattice model then often there
are dilute versions of the model which are also solvable.
These dilute models are studied in \cite{MR94c:82032},
\cite{MR94m:20084}, \cite{MR95i:82031} and \cite{MR94f:82024}.
In this paper we discuss the algebras generated by the single bond
transfer matrices and their representation theory.

Our motivation for studying these dilute models is that in
\cite{MR1984741} it is shown that there is a connection between
the dilute Potts models and the exceptional series of Lie
algebras proposed in \cite{MR96m:22012}, \cite{MR97k:22008}
and \cite{MR97i:17005}.

\section{Categories}
In this section we give the main construction of this paper.
This construction gives a tensor product of monoidal categories.
Before giving this construction we make some remarks on monoidal
categories.

In this paper nearly all the monoidal categories are
$K$-linear for some integral domain $K$ and furthermore these two
structures are compatible. In practice, it is always the case  that
$K=\End(1)$.

If $\cC$ is a category then we can 
construct $K\cC$, the free $K$-linear category on $\cC$.
The universal property of this construction is that it gives
a left adjoint to the forgetful functor from $K$-linear categories
to categories. This construction passes to monoidal categories.

Another useful property of a $\Z$-linear category is the additive
property. This property says that we can form the direct sum of two
(or more) objects and that there is a zero object.
If $\cC$ is a $\Z$-linear category then we can construct a free
additive category on $\cC$ by taking objects to be vectors of
objects in $\cC$ and morphisms to be matrices of morphisms.
The universal property of this construction is that it gives
a left adjoint to the forgetful functor from additive categories
to $\Z$-linear categories. This construction passes to monoidal categories.

The following construction defines a tensor product of 
categories. This construction gives a category which is equivalent
to the construction in \cite[1.1.15]{MR2002d:18003}.
Both of these constructions are tensor products in the sense of
\cite[\S 5]{MR92d:14002}.
\begin{definition}\label{tens}
Let $\cC_1$ and $\cC_2$ be two monoidal categories. Then we construct
a new monoidal category, $\cC_1\otimes \cC_2$ as follows. The objects
are sequences of elements of the disjoint union $\Ob(\cC_1)\cup\Ob(\cC_2)$.
Let $X$ be any such sequence. Then associated to $X$ is a sequence
$X_1$ of elements of $\Ob(\cC_1)$ and a sequence $X_2$ of elements of
$\Ob(\cC_2)$. Then given two sequences $X$ and $Y$ we define a
morphism $\phi\colon X\rightarrow Y$ to be a pair $(\phi_1,\phi_2)$
where $\phi_1\colon \otimes_{x\in X_1} x \rightarrow \otimes_{y\in Y_1} y$
is a morphism in $\cC_1$ and similarly
$\phi_2\colon \otimes_{x\in X_2} x \rightarrow \otimes_{y\in Y_2} y$
is a morphism in $\cC_2$. The tensor product is defined on objects
by concatenating sequences and is defined on morphisms by taking the
tensor products in $\cC_1$ and $\cC_2$.
\end{definition}

If $\cC_1$ and $\cC_2$ are both $K$-linear monoidal categories then
we modify this construction so that it gives another $K$-linear 
category. In this version we take
\[
\Hom(X,Y) = \Hom_{\cC_1}(\otimes_{x\in X_1}x,\otimes_{y\in Y_1}y)
\otimes \Hom_{\cC_2}(\otimes_{x\in X_2}x,\otimes_{y\in Y_2}y)
\]

If $\cC_1$ and $\cC_2$ are both $K$-linear and additive monoidal categories
then we take the free additive category on this to give a monoidal
category which is also $K$-linear and additive.

The first property of this construction is that it is associative.
This means that we natural equivalences of monoidal categories
\begin{equation}
(\cC_1\otimes \cC_2)\otimes \cC_3 \sim
\cC_1\otimes (\cC_2\otimes \cC_3)
\end{equation}
The simplest way to see this is to note that both sides are naturally
equivalent to a monoidal category whose objects are sequences of
elements of 
$\Ob(\cC_1)\cup\Ob(\cC_2)\cup\Ob(\cC_3)$.

A further property is that it is symmetric. This structure is
given by natural equivalences of categories
\begin{equation}\label{symm}
\cC_1 \otimes \cC_2 \cong \cC_2 \otimes \cC_1
\end{equation}
This functor is given by
\[ (\phi_1,\phi_2) \mapsto (\phi_2,\phi_1) \]

Note that if $H_1$ and $H_2$ are Hopf algebra and $\Rep(H_1)$
and $\Rep(H_2)$ are monoidal categories of representations then
we have a natural functor
\[ \Rep(H_1) \otimes \Rep(H_2) \rightarrow \Rep(H_1\otimes H_2) \]
This functor is defined on objects as follows. Each object $V$
of $\Rep(H_1)$ gives an object $V\otimes 1$ of $\Rep(H_1\otimes H_2)$
and each object $V$ of $\Rep(H_2)$ gives an object $1\otimes V$ of
$\Rep(H_1\otimes H_2)$. Then each object of $\Rep(H_1) \otimes \Rep(H_2)$
can be regarded as a sequence of objects of $\Rep(H_1\otimes H_2)$
and we take the tensor product to give an object of $\Rep(H_1\otimes H_2)$.

Here we recall some standard constructions in a monoidal category.
\begin{definition}
Let $\cC$ be a $K$-linear monoidal category and $V$ an object of $\cC$.
Then $\cI(V)$ is the category with objects $\{n|n\ge 0\}$ and morphisms
given by
\[ \Hom_{\cI(V)}(n,m) = \Hom_\cC(\otimes^nV,\otimes^mV) \]
\end{definition}
This is a $K$-linear monoidal subcategory of $\cC$ and is called the
category of invariant tensors.
Closely related to this is the sequence of algebras
$\{A(n)\}$ given by
\begin{equation}\label{cons}
A(n) = \End_{\cC}(\otimes^nV) 
\end{equation}
which are the endomorphism algebras of the objects of the category
of invariant tensors.
This construction is the motivation for the
following definition:

\begin{definition}\label{tow}
A tower of algebras is a sequence of algebras
$\left\{ A(n)\right\} $ together with homomorphisms
$\phi_{n,m}:A(n)\otimes A(m)\rightarrow A(n+m)$ which satisfy the
associativity condition.
\[ \varphi_{r,s+t}(1\otimes\varphi_{s,t}) =
\varphi_{r+s,t}(\varphi_{r,s}\otimes1) \]
\end{definition}
Another way of stating this associativity condition is to say
that the following diagram commutes:
\begin{equation*}
\begin{CD}
A(r)\otimes A(s)\otimes A(t) @>{1\otimes\varphi_{s,t}}>> A(r)\otimes A(s+t) \\
@V{\varphi_{r,s}\otimes 1}VV @VV{\varphi_{r,s+t}}V \\
A(r+s)\otimes A(t) @>>{\varphi_{r+s,t}}> A(r+s+t) \\
\end{CD}
\end{equation*}
Then it is clear that the sequence of algebras in (\ref{cons}) is a
tower of algebras; this uses the tensor product in $\cC$.

We can also regard a tower of algebras as a $K$-linear monoidal category
denoted by $\cA$.
\begin{definition}\label{cat}
The objects of $\cA$ are $\{n|n\ge 0\}$ and morphisms
are given by
\begin{equation} \Hom_{\cA}(n,m) = 
\left\{\begin{array}{cl}
A(n) & \text{if $n=m$} \\
0 & \text{otherwise}
\end{array}\right. \end{equation}
\end{definition}

Let $V$ be an object of a monoidal category $\cC$. Then we can
take the tower of algebras associated to $V$ in (\ref{cons})
and then construct the category $\cA(V)$ by (\ref{cat}).
Note that we then have an inclusion of $K$-linear monoidal
categories $\cA(V) \rightarrow \cI(V)$.

For applications to knot theory and topological field theory
these categories are required to have extra structure. Next
we show that these structures pass to the tensor product.

Examples of tensor categories are tensor categories
defined by diagrams. The main examples are the Temperley-Lieb
category in \cite{MR96h:20029}, the category of braids, and the
various categories of tangles (tangles, oriented tangles,
framed tangles, and oriented framed tangles). Let $\cD$ be
one of these tensor categories of diagrams. Then consider the
repeated tensor product $\otimes^c\cD$. This can be considered 
as a diagram category where a diagram consists of a diagram
with strings labelled by elements of $C$. Then taking the
strings labelled by $c\in C$ gives a diagram in $\cD$ and
strings of different colours can cross. This is illustrated
in \cite{math-ph/0307017}. Then it is observed in \cite{MR94m:20084}
that there are functors of tensor categories
\begin{equation}\label{braid}
\cD \rightarrow \otimes^c (\Z\cD)
\end{equation}
for all $c$. Let $D$ be a diagram in $\cD$. Then the result of applying
the functor to $D$ is the sum of all possible colourings of $D$.

The Temperley-Lieb category depends on a parameter whereas the
other examples do not involve a parameter. Next we explain how
this parameter behaves when this construction is applied to
Temperley-Lieb categories.
Let $\cT(\delta)$ be the Temperley-Lieb category with parameter
$\delta$. Then note that this construction gives functors
\[
\cT(\sum_{c\in C}\delta_c) \rightarrow \otimes_{c\in C} \cT(\delta_c)
\]

\section{Basic Example}
Next we take a very simple monoidal category and apply the
construction in Definition (\ref{tens}) to several copies
of this category.

\begin{definition}
Let $\cF$ be the category with objects $n$ and morphisms given by
\[ \Hom_{\cC}(n,m) =
\left\{\begin{array}{cl}
\{1\} & \text{if $n=m$} \\
0  & \text{otherwise}
\end{array}\right. \]
This is a monoidal category with tensor product given by
$n\otimes m=n+m$.
\end{definition}

The main property of this example is that the object $1$ in the
monoidal category $\cF$ is a universal object in a monoidal 
category. This means that if $V$ is an object in a monoidal
category $\cC$ then there is a unique monoidal functor 
$\cF\rightarrow\cC$ which sends the object $1$ to $V$.
A consequence of this and the naturality of the product in
Definition (\ref{tens}) is that if we have objects $V_1$ and
$V_2$ in monoidal categories $\cC_1$ and $\cC_2$ then we have
a monoidal functor 
\[ \cF \otimes \cF \rightarrow \cC_1 \otimes \cC_2 \]

Let $\Z\cF$ be the free $\Z$-linear category on $\cF$.
Let $\rho$ be a one dimensional representation of $\Z$ such that
$\otimes^n\rho \ne 1$ for all $n>0$.
Then $\Z\cF$ is also the category of invariant tensors
associated to the representation $\rho$ of $\Z$.

Then we consider the repeated product $\otimes^c\cF$.
This category is a groupoid. Let $C$ be a set with $c$ elements.
The objects are sequences of elements
of $C$. Given a sequence $X$, let $p\colon C\rightarrow \N$ be the
function such that $p(c)$ is the number of times $c$ appears in
the sequence $X$. Then two sequences are isomorphic if and only
if they have the same function. Also if two objects are isomorphic
then there is only one isomorphism. One way to see this is to take
the morphisms to be permutation diagrams with edges labelled
by $C$ and such that if two edges cross then they are labelled by
distinct colours.

Now consider the repeated product $\otimes^c(\Z\cF)$ and take
the free additive category on this repeated tensor product.
Recall that the objects of the repeated product $\otimes^c(\Z\cF)$ are
sequences of integers each labelled by a colour.
\begin{definition}\label{bas} Let $V$ be the sum of the
$c$ sequences of length one consisting of the single integer 1.
Consider the category of invariant tensors
of $V$ and let $\{F^{(c)}(n)\}$ be the associated tower of algebras.
\end{definition}
Given a function
$p\colon C\rightarrow \N$ let $|p|=\sum_{c\in C}p(c)$.

\begin{lemma}\label{repo}
For all $c\ge 1$ and $n\ge 0$ the algebra $F^{(c)}(n)$ is a direct
sum of matrix algebras. The simple
$F^{(c)}(n)$-modules are indexed by functions, $p$, such that
$|p|=n$. The dimension of the simple module associated to $p$
is the multinomial coefficient $\binom{|p|}{\{p(c)|c\in C\}}$.
\end{lemma}

\section{Bratteli diagrams}
Let $V$ be an object in a monoidal category $\cC$, then associated
to $V$ is the category of invariant tensors and the tower of algebras.
In this section we let $V_1$ be an object in $\cC_1$ and 
$V_2$ an object in $\cC_2$.
\begin{definition}\label{dou}
Consider both $V_1$ and $V_2$
as objects in $\cC_1\otimes \cC_2$ by taking them to be sequences
of length one. Assume $\cC_1\otimes \cC_2$ is additive 
and put $V=V_1\oplus V_2$.
\end{definition}
Then we consider the category of invariant tensors associated to
$V$ and the tower of algebras associated to $V$.

First we compare $\cI(V)$ with $\cI(V_1)\otimes\cI(V_2)$.
The objects of the tensor
product $\cI(V_1)\otimes \cI(V_2)$ are sequences where each term is
of the form $\otimes^nV_1$ or $\otimes^nV_2$ for some $n\ge 0$.
Let $\cI$ be the full subcategory with objects those sequences
in which every term is of the form $V_1$ or $V_2$. Then by
construction $\cI$ is a full monoidal subcategory of
$\cI(V_1)\otimes \cI(V_2)$. The connection between $\cI$ and $\cI(V)$
comes from noting that $\otimes^nV$ is the direct sum of
of all objects of $\cI$ which are sequences of length $n$.

Next we note that the inclusion functor $\cI\rightarrow 
\cI(V_1)\otimes \cI(V_2)$ is an equivalence of categories.
The inverse functor is constructed
by writing a term $\otimes^nV_1$ by $n$ terms $V_1$ and similarly
for $\otimes^nV_2$.

The object $V$ can be considered as an object of 
$\cC_1\otimes\cC_2$, $\cI(V_1)\otimes \cI(V_2)$ or $\cI$.
All these three cases give the same category of invariant tensors
and the same tower of algebras. Let $A(n)$ be the endomorphism algebra
of $\otimes^nV$ cosidered as an object in any of these categories.
Another possibility is to
consider $V$ as an object of $\cA(V_1)\otimes \cA(V_2)$.
\begin{definition}
Define $\widehat{A}(n)$ to be the endomorphism algebra
of $\otimes^n V$ considered as an object of $\cA(V_1)\otimes \cA(V_2)$.
\end{definition}
Then, by construction, we have inclusions $\widehat{A}(n)
\rightarrow A(n)$ for all $n\ge 0$.
If the inclusions $\cA(V_1)\rightarrow \cI(V_1)$ and
$\cA(V_2)\rightarrow \cI(V_2)$ are both isomorphisms then
these inclusions will be isomorphisms. However, in general,
$\widehat{A}(n)$ will be a proper subalgebra of $A(n)$.

Here we consider this tower of subalgebras. The main result is:
\begin{proposition}\label{dah} For $n\ge 0$, the algebra $\widehat{A}(n)$
is isomorphic to 
\begin{equation*}\label{ts}
\bigoplus_{r+s=n} M(\binom{n}{r,s})\otimes A_1(r)\otimes A_2(s)
\end{equation*}
where $M(N)$ is the algebra of $N\times N$ matrices.
\end{proposition}
\begin{proof}
First note that the inclusion $F^{(2)}(n)\rightarrow \widehat{A}(n)$
is an inclusion
\[
\bigoplus_{r+s=n} M(\binom{n}{r,s}) \rightarrow \widehat{A}(n)
\]
Now let $e$ be a diagonal elementary matrix in $M(\binom{n}{r,s})$,
so $e$ is an idempotent permutation diagram.
Then observe that we have an isomorphism
\[ e\widehat{A}(n)e \cong A_1(r)\otimes A_2(s) \]
These idempotents give a decomposition of the identity on both
sides of (\ref{ts}) so the result follows since both of these
decompositions give the same Pierce decomposition.
\end{proof}

In this section we consider the Bratteli diagrams of the towers
of algebras we are considering. The Bratteli diagram was introduced
in \cite{MR47:844} in the study of approximately finite $C^*$-algebras.
Assume we are given a tower of algebras such that each
algebra is a direct sum of matrix algebras.
Then associated to this tower of algebras is a graded directed graph
called the Bratteli diagram. This has the following properties.
The vertices of degree $n$ correspond to the (isomorphism classes of)
simple $A_1(n)$-modules; so, there is a single vertex $v_0$ of degree
$0$. Let $v$ be a vertex of degree $n$ associated to the simple
module $M$. Then the dimension of $M$ is the number of directed paths
from $v_0$ to $v$.

\begin{definition}\label{teng}
Next let $\Gamma_1$ be a directed graph with edge set $E_1$, vertex
set $V_1$ and maps $h_1,t_1\colon E_1\rightarrow V_1$. Let $\Gamma_2$
be second directed graph. Then we define the product to have edge
set $E=(V_1\times E_2)\cup (E_1\times V_2)$, vertex set $V=V_1\times V_2$
and define $h,t\colon E\rightarrow V$ by
\begin{equation*}
\begin{array}{cc}
h(e,v)=(h_1(e),v) & h(v,e)=(v,h_2(v)) \\
t(e,v)=(t_1(e),v) & t(v,e)=(v,t_2(v)) 
\end{array}
\end{equation*}
\end{definition}

Note that this product is also associative and symmetric.
Furthermore if $\Gamma_1$ and $\Gamma_2$ have basepoints $v_1$
and $v_2$ then we take $(v_1,v_2)$ to be the basepoint in
$\Gamma$.

Assume that the two towers $\{A_1(n)\}$ and $\{A_2(n)\}$ both have
Bratteli diagrams. Then the tower of algebras $\{\widehat{A}(n)\}$ also has
a Bratteli diagram and this Bratteli diagram is the product of the
two Bratteli diagrams given in Definition (\ref{teng}).

This observation and Proposition (\ref{dah}) both give
the following dimension formula. An irreducible
representation, $W$, of $\widehat{A}(n)$ is labelled by a pair $(W_1,W_2)$
where $W_1$ is an irreducible representation of $A_1(r)$ and
$W_2$ is an irreducible representation of $A_2(s)$ where
$r+s=n$. Then the dimension of $W$ is given by
\begin{equation}\label{dimp}
\dim(W) = \binom{n}{r,s} \dim(W_1)\dim(W_2)
\end{equation}

There are two examples in which the inclusions $\widehat{A}(n)
\rightarrow A(n)$ are isomorphisms. The first is the basic example
of the object 1 in $\Z\cF$. Then this construction shows that the
Bratteli diagram for the tower of algebras $\{F^{(2)}(n)\}$ is
Pascal's triangle. More generally, the Bratteli diagram for the
tower of algebras $\{F^{(c)}(n)\}$ is the generalisation of Pascal's
triangle which gives the multinomial coefficients. This directed
graph has vertices $\N^c$ with directed edges given by increasing
a single coordinate by 1 and the basepoint is the origin.

For the second example we interpret the Hecke algebras as endomorphism
algebras of tensor powers of some object in a monoidal category.
This can be achieved as follows. Let $R$ be the ring defined by
\[ R = \Z[\delta ,q ,z, 1/qz ]/
\langle \delta(q-q^{-1})=(z-z^{-1}) \rangle \]
Then let $\cH$ be the $R$-linear monoidal category obtained by
taking the free $R$-linear category on the category of oriented
tangles and imposing the HOMFLY skein relations, see
\cite{MR89c:46092}. Then let $V$ be the object whose identity morphism is
a single descending string. Then the Hecke algebra $H(n)$ is the
endomorphism algebra of $\otimes^n V$. Then in this example we also
have that the inclusion $\widehat{H}(n) \rightarrow H(n)$ is an
isomorphism. Then this construction shows that the Bratteli diagram
for the tower of algebras $\{H^{(c)}(n)\}$ is the Bratteli diagram
for the Ariki-Koike algebras given in \cite{MR95h:20006}.

\section{Symmetry}\label{sym}
Let $V$ be an object of $\cC$. Then $V$ can be regarded as an
object in $\otimes^c\cC$ in $c$ different ways and we let
$V^{(c)}$ be the direct sum of these. For $c=2$ this object
is given by taking $V_1=V_2$ in Definition (\ref{dou}).
Let $\{A^{(c)}(n)\}$ be the tower of algebras associated to $V^{(c)}$.
For example, if we take $V$ to be the object 1 of $\Z\cF$ then
this construction gives the algebras in Definition (\ref{bas}).
Then the tower of algebras has an action of the symmetric group
$S_c$ which arises from the symmetric structure on the product
in Definition (\ref{tens}). 
Let $A^{S(c)}(n)$ be the tower of algebras obtained by
taking the fixed point subsets. Then the symmetric group acts freely
so we have
\begin{equation}\label{dimi}
\dim A^{S(c)}(n) = \frac1{c!} \dim A^{(c)}(n)
\end{equation}
for all $n>0$.

Assume $V$ is a representation of a Hopf algebra $H$. Then
$V^{(c)}$ is the representation of $\otimes^cH$ given by
\[
\bigoplus_{i=1}^c
\overbrace{1\otimes \ldots \otimes 1}^{\text{$i-1$ factors}}
\otimes V \otimes
\overbrace{1\otimes \ldots \otimes 1}^{\text{$c-i$ factors}}
\]
Then, for $n\ge 0$, the algebra $A^{(c)}(n)$ is the endomorphism
algebra of $\otimes^nV^{(c)}$.
The tower of algebras $\{A^{S(c)}(n)\}$ also has an interpretation
as the tower associated to a representation of a Hopf algebra.
This interpretation is given by using the same representation $V^{(c)}$
and regarding it as a representation of the wreath product
$S_c \wr H$.

Next we give a generalisation of the construction (\ref{braid}).
Let $\cC$ be a $K$-linear monoidal
category and assume $\cC\otimes \cC$ is additive. 
Then there is a functor $\cC \rightarrow \cC\otimes \cC$ which
on objects is given by $V \mapsto V^{(2)}$ and which on
morphisms is given by $\phi \mapsto \phi \oplus \phi$.

In particular we have algebra homomorphisms
\[ A(n) \rightarrow A^{(c)}(n) \]
These factor through $A^{S(c)}(n)$ to give algebra homomorphisms
\[ A(n) \rightarrow A^{S(c)}(n) \]

Next we discuss the representation theory of the algebras
$\{A^{S(2)}(n)\}$. Let $s$ be the involution of $A^{(2)}(n)$
whose fixed point set is $\{A^{S(2)}(n)\}$. Then $s$ induces
an involution on the representations of $\{A^{(2)}(n)\}$.
Let $W$ be an irreducible representation of $\{A^{S(2)}(n)\}$.
Then if $W$ and $s(W)$ are not isomorphic $W$ can be regarded as  
an irreducible representation of $\{A^{S(2)}(n)\}$; and if $W$
and $s(W)$ are isomorphic $W$ is the direct sum of two 
irreducible representations of $\{A^{S(2)}(n)\}$ of equal dimension.

Now assume that $\{A^{(2)}(n)\}$ is a direct sum of matrix algebras
so that the dimension of $\{A^{(2)}(n)\}$ is the sum of the squares
of the dimensions of the irreducible representations. Then 
$\{A^{(2)}(n)\}$ is also a direct sum of matrix algebras and the
sum of the squares of the dimensions is half the dimension 
of $\{A^{(2)}(n)\}$ which is consistent with (\ref{dimi}).

\begin{center}
\begin{figure}
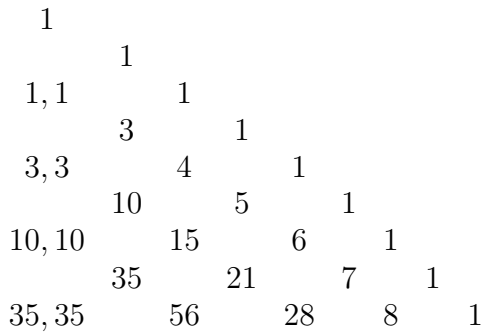

\[ \begin{array}{ccccccccc}
1 & & &  & &  & &  & \\
 & 1 & &  & &  & &  & \\
1,1 & & 1 & &  & &  & & \\
 & 3 & & 1 & &  & &  & \\
3,3 & & 4 & & 1 & &  & & \\
 & 10 & & 5 & & 1 & &  & \\
10,10 & & 15 & & 6 & & 1 & & \\
 & 35 & & 21 & & 7 & & 1 & \\
35,35 & & 56 & & 28 & & 8 & & 1
\end{array} \]
\caption{Bratteli diagram of $\{F^{S(2)}(n)\}$.}\label{bra}
\end{figure}
\end{center}

As an example of this we consider the tower of algebras $F^{S(2)}(n)$.
The Bratteli diagram for these algebras is given in Figure \ref{bra}.
This is essentially Pascal's triangle folded about the central axis.
This is also the Bratteli diagram for the Temperley-Lieb algebras
of type $D_n$ (see \cite{MR2000i:57006} and \cite{MR97k:20065}).
Each two part partition of $n$ can be written as $(n-p,p)$ where
$0\le 2p\le n$. Corresponding to each such partition there is a
representation of the $n$-string algebra whose dimension is given
by the binomial coefficient $\binom{n}{p}$. These representations
are irreducible except if $2p=n$ in which case they are the sum of
irreducible representations of the same dimension.

\begin{definition}
Let $H$ be the group algebra of the wreath product $ \Z \wr S_2 $.
This is a Hopf algebra. The algebra is obtained from the ring of
Laurent polynomials $K[q_1,q_2,1/q_1q_2]$ by adjoining an element
$\sigma$ which satisfies the relations
\[ \sigma^2=1 \qquad \sigma q_1 = q_2\sigma \qquad \sigma q_2
= q_1 \sigma \]
The coproduct is given by
\[
\Delta(q_1)=q_1\otimes q_1 \qquad
\Delta(q_2)=q_2\otimes q_2 \qquad
\Delta(\sigma)=\sigma\otimes\sigma 
\]
and the antipode is given by
\[ S(q_1)=q_1^{-1} \qquad S(q_2)=q_2^{-1} \qquad S(\sigma)=\sigma \]
\end{definition}

Note that we can eliminate the generator $q_2$ using 
$q_2=\sigma q_1\sigma$. Then the defining relations become
$\sigma^2=1$ and $\sigma q_1\sigma q_1 = q_1\sigma q_1\sigma$.

\begin{definition}
Let $x$ be an invertible scalar.
Then we define a two dimensional representation $\rho^{(n)}$ for
$n\ge 0$ by
\begin{equation}\label{rep2}
q_1 \mapsto 
\left(\begin{array}{cc}
x^n & 0 \\
0 & 1
\end{array}\right)
q_2 \mapsto 
\left(\begin{array}{cc}
1 & 0 \\
0 & x^n
\end{array}\right)
\sigma \mapsto 
\left(\begin{array}{cc}
0 & 1 \\
1 & 0
\end{array}\right)
\end{equation}
\end{definition}

Then $\rho^{(n)}$ is irreducible if $x^n\ne 1$ and $\rho^{(0)}$
is the sum of two one dimensional representations. Denote these
by $\rho_+$ and $\rho_-$. The tensor
product of these representations with $\rho^{(1)}$ are given by the
following two Lemmas.
\begin{lemma} For all $n>0$,
\[
\rho^{(1)} \otimes \rho^{(n)} = \rho^{(n-1)} \oplus \rho^{(n+1)}
\]
\end{lemma}
\begin{proof} Let the ordered basis of $\rho^{(n)}$ used in
(\ref{rep2}) be $(e_1^{(n)},e_2^{(n)})$. 
Consider the tensor
product representation. The action on the subspace with
ordered basis $(e_1^{(1)}\otimes e_1^{(n)},e_2^{(1)}\otimes e_2^{(n)})$
is the representation $\rho^{(n+1)}$ and the action on the
subspace with ordered basis
$(\frac1x e_2^{(1)}\otimes e_1^{(n)},\frac1x e_1^{(1)}\otimes e_2^{(n)})$
is the representation $\rho^{(n-1)}$.
\end{proof}
\begin{lemma}
\[ \rho^{(1)} \otimes \rho_{\pm} = \rho^{(1)} \]
\end{lemma}
\begin{proof} The representation $\rho_+$ is the trivial
representation so this case is clear. The representation
$\rho^{(1)}\otimes\rho_-$ is given by taking the matrices in
(\ref{rep2}) for $n=1$ and multiplying the matrix representing
$\sigma$ by $-1$.
Let the ordered basis of this representation
be $(f_1,f_2)$. Then if we change to the ordered basis
$(f_1,-f_2)$ we get the representation $\rho^{(1)}$.
\end{proof}
These tensor product decompositions show that if $x^n\ne 1$ for
$n>0$ then the tower of
algebras given by $\End_{H}(\otimes^n\rho^{(1)})$ has the Bratteli
diagram given in Figure \ref{bra}.

\section{Dilute Temperley-Lieb}
In this section we consider the $c$-colour Temperley-Lieb algebras,
$T^{(c)}(n)$.
Then it is clear from the diagram point
of view that each $T^{(c)}(n)$ is a cellular algebra in the sense of
\cite{MR97h:20016}. Next we discuss this cell structure in more detail,
following \cite{math-ph/0307017}.
Each diagram $D\in T^{(c)}(n)$ has say $|p|$ propagating strings. Then
these $|p|$ strings are coloured, so, say $p(c)$ are coloured by
$c\in C$. Then $p$ is a function $p\colon C\rightarrow \N$.
Define a partial order on these functions by $p_1\le p_2$
if $p_1(c)\le p_2(c)$ for all $c\in C$. Then let $I(p)$ be the
subspace with basis all diagrams, $D$, such that $p(D)\le p$.
Then for each $p$, $I(p)$ is an ideal and $I(p_1)\subset I(p_2)$
if and only if $p_1\le p_2$.

The algebra $F^{(c)}(n)$ is also a quotient of $T^{(c)}(n)$.
The tower of algebras, $\{T^{(c)}(n)\}$, is obtained from the tower 
$\{F^{(c)}(n)\}$ by a Jones tower construction.
This construction is described in
\cite[Chapter 2]{MR91c:46082} and \cite[\S 4]{math.RT/0401314}.
The Bratteli diagram for the tower of algebras $F^{(c)}(n)$
is a directed graph with vertices $\N^c$. Then we obtain the
Bratteli diagram for the tower of algebras $T^{(c)}(n)$ by
taking paths in the underlying undirected graph. In particular,
the irreducible representations of $T^{(c)}(n)$ are indexed by
sequences $(k_1,\ldots ,k_c)$ such that $k_i\ge 0$ and
$n-k_1- \cdots k_c$ is even and non-negative.
For the two-colour case, these
paths are counted in \cite{MR93i:05008} and \cite{MR2001d:05006}.

The dimensions of the two colour algebras are given by
\[ \begin{array}{cccc}
0 & 1 & 2 & 3  \\
1 & 2 & 10 & 70
\end{array} \]
The simplest description of these numbers is the formula
\begin{equation}\label{dim}
\dim T^{(2)}(n) = C(n)C(n+1) 
\end{equation}
where $C(n)$ is the Catalan number.
For $c$ colours we have the formula
\[ \dim T^{(c)}(n) = \sum_{n_1+\cdots n_c = n}
\left( \begin{array}{c} 2n \\ 2n_1,\ldots ,2n_c
\end{array} \right)
\prod_{i=1}^c C(n_i) \]
Equivalently, in terms of exponential generating functions, 
\[ \sum_{n\ge 0} \frac{\dim T^{(c)}(n)}{(2n)!}z^{2n} =
 \left( \sum_{n\ge 0} \frac{C(n)}{(2n)!}z^{2n} \right)^c
 \]

This can be extended to the dimensions of the irreducible
representations. let $F(x,y)$ be defined by
\[ F(x,y) = \sum_{n,p\ge 0}
\frac{(n-2p+1)}{p!(n-p+1)!}x^ny^{n-2p} \]
Let $S(n;k_1,\ldots ,k_c)$ be the simple $T^{(c)}(n)$-module
associated to the vector $(k_1,\ldots ,k_c)$ and define 
\[ G(x;y_1,\ldots ,y_c) =
\sum_{n,k_1,\ldots ,k_c}
\dim S(n,k_1,\ldots ,k_c)\frac{x^n}{n!}y_1^{k_1}\ldots
y_c^{k_c}
\]
Then these are related by
\[ G(x;y_1,\ldots ,y_c) = \prod_{i=1}^c F(x,y_i) \]

In this rest of this section we restrict attention
to just two colours and we also assume that the two parameters
$\{ \delta(c) | c \in C\}$ are both equal and we denote them
both by $\delta$. Then we let $T^{S(2)}(n)$ be the fixed point
subalgebra of the involution which interchanges the two colours.

The general discussion in \S \ref{sym} applies and gives the following
description of the irreducible representations of $T^{S(2)}(n)$.
Let $S(n;r,s)$ be an irreducible representation of $T^{(2)}(n)$.
Then if $r\ne s$ the restriction to $T^{S(2)}(n)$ is irreducible
and the restrictions of $S(n;r,s)$ and $S(n;s,r)$ are isomorphic.
The restriction of each representation $S(n,r,r)$ is the direct
sum of two irreducible representations of the same dimension.

\begin{definition} Define the following elements of $F^{(2)}(n)$.
The element $\pi_i(c)$ is the sum over all colourings of the
identity permutation with the string $i$ coloured by $c$.
The element $s_i(c,d)$ is the sum over all colourings of
the permutation $(i,i+1)$ such that the string from $i$ to
$i+1$ is coloured $c$ and the string from $i+1$ to $i$ is coloured
$d$.
\end{definition}

\begin{definition} Define the following elements of $T^{(2)}(n)$.
Let $U_i$ be the diagram for the standard generators of the
Temperley-Lieb algebras. Then the element $u_i(c,d)$ is obtained
by summing over all colourings such that the top arc is coloured
$c$ and the lower arc is coloured $d$.
\end{definition}

\begin{definition}\label{fix}
For $n>1$, define the following elements of $T^{S(2)}(n)$,
\begin{gather*}
e_i = \pi_i(1)\pi_i(1) + \pi_i(2)\pi_i(2) \\
f_i = \pi_i(1)\pi_i(2) + \pi_i(2)\pi_i(1) \\
s_i = \sigma_i(1,2) + \sigma_i(2,1) \\
u_i = u_i(1,1)+u_i(2,2) \\
t_i = u_i(1,2)+u_i(2,1) 
\end{gather*}
\end{definition}
Then these elements generate $T^{S(2)}(n)$. For fixed $i$, the
five dimensional algebra with these elements as basis also has a
basis given by the following five orthogonal idempotents:
\begin{equation}\label{oidem}
e_i-\frac1\delta u_i \qquad
\frac12 (f_i\pm s_i) \qquad
\frac1{2\delta}(u_i\pm t_i)
\end{equation}

In order to describe the $R$-matrices and the braid matrices we
extend the ring of scalars from the polynomial ring $\Z[\delta]$
to the ring of Laurent polynomials $\Z[q,1/q]$ using the ring
homomorphism determined by
\begin{equation}\label{spec}
\delta \mapsto -q^2-q^{-2}
\end{equation}

The construction (\ref{braid}) shows that the algebras
$\{T^{S(2)}(n)\}$ can be used to construct an invariant
of unoriented framed links just as the Temperley-Lieb
algebras can be used to construct the Kauffman bracket
polynomial. This invariant has the following description.
First, for a link $L$ let $\langle L\rangle$ be the Kauffman
bracket normalised so that the empty link has the invariant 1.
Then this invariant is multiplicative under disjoint union.
Then the two-colour dilute version of this invariant is
\[ \sum_{L=L_1\cup L_2} \langle L_1\rangle\langle L_2\rangle \]
The sum is over all ways of colouring the components by the two
colours 1 and 2; and for a given colouring $L_1$ is the sublink
coloured 1 and $L_2$ is the sublink coloured 2. Although they
are different invariants, there are some similiarities with
the link invariants in \cite{MR92k:57011}. Next we show that
this invariant can be calculated by taking a Markov trace on the
algebras $\{T^{S(2)}(n)\}$.

\begin{proposition}
The tower of algebras $\{T^{S(2)}(n)\}$ has the property that
\begin{multline}\label{bim}
T^{S(2)}(n+1) = 
T^{S(2)}(n) + T^{S(2)}(n)e_nT^{S(2)}(n) \\
+ T^{S(2)}(n)s_nT^{S(2)}(n)
+ T^{S(2)}(n)u_nT^{S(2)}(n) + T^{S(2)}(n)t_nT^{S(2)}(n)
\end{multline}
\end{proposition}
\begin{proof} First we check that the necessary condition in
\cite{MR98f:57018} is satisfied. Then this can be shown by writing
down a long, but finite, list of relations.
\end{proof}

The map of the braid group given in (\ref{braid}) is
\begin{equation}\label{hom}
\sigma_i^{\pm 1} \mapsto q^{\pm 2}e_i + u_i - Q^{\pm 1}s_i 
\end{equation}
where $Q$ is an independent parameter. This parameter arises
by taking into account the fact that the defining relations
of the braid group are homogeneous so any representation can
be scaled.

The map of the Temperley-Lieb algebra $T(n)\rightarrow
T^{S(2)}(n)$ is given by
\[ U_i \mapsto u_i+t_i \]
Note that the value of a closed loop in $T(n)$ is $2\delta$.
There is a conditional expectation
$\varepsilon_n\colon T^{S(2)}(n+1)\rightarrow T^{S(2)}(n)$
This conditional expectation is determined by
\[ U_{n+1}aU_{n+1} = \varepsilon_n(a)U_{n+1} \]
It follows from (\ref{bim}) that this is determined by
\begin{equation*}
\varepsilon_n(a) = \delta a,
\varepsilon_n(ae_nb) = \delta ab,
\varepsilon_n(as_nb) = 0,
\varepsilon_n(au_nb) = ab,
\varepsilon_n(at_nb) = 0
\end{equation*}
for all $a,b\in T(n)$.
This give a sequences of traces $\tau_n\colon T(n)\rightarrow \Z[\delta]$
which are determined by $\tau_{n+1}(a) = \tau_n(\varepsilon_n(a))$
for all $a\in T(n+1)$. This sequence of traces gives a sequence of traces
of the braid group algebras which satisfies the Markov property
since we have
\[ \varepsilon_n(a\sigma_n^{\pm 1}b) = (q^{\pm 2}\delta +1)ab
= -q^{\pm 4}ab \]
using (\ref{spec}).

The Yang-Baxter equation is the equation
\begin{equation}\label{YB}
R_i(u)R_{i+1}(uv)R_i(v)=R_{i+1}(v)R_i(uv)R_{i+1}(u)
\end{equation}
for $|i-j|>1$ and all $u$ and $v$.
The solution to the Yang-Baxter equation is given in
\cite{MR94f:82024} and \cite[(4.34)]{MR94m:20084}.

Introduce the notation
\[
[ax+b] =\frac{ q^{b}u^{a} -q^{-b}u^{-a} }{q-q^{-1} }
\]
This $R$-matrix is given in terms of the elements in
Definition (\ref{fix}) by:
\begin{multline}
R_i(u) = 
[1-x][3-x]e_i
+[3-x]f_i \\
-[x][2-x]u_i
+[x]t_i
+[x][3-x]s_i
\end{multline}
The relation (\ref{YB}) can be checked by checking it each
irreducible representation of $T^{S(2)}(3)$. These have dimension
1, 3 and 5.

This solution has a number of properties. The first is that
\[ R_i(1) = [3] \]
Another property is that
taking the coefficients of $u^{\pm 2}$ gives a representation of the
braid group. These coefficients are given by 
\[ \frac{q^{\mp 2}}{(q-q^{-1})^2}\sigma_i^{\mp 1} \]
where $\sigma_i^{\pm 1}$ is given by (\ref{hom}) with $Q=q$.

Another property of this $R$-matrix is that it has crossing symmetry.
This means that it is invariant under the involution
\begin{equation*}
x \leftrightarrow 3-x \quad
u_i \leftrightarrow e_i \quad
t_i \leftrightarrow f_i \quad
s_i \leftrightarrow s_i 
\end{equation*}
The subalgebra generated by the braid group is the subalgebra
generated by $\{e_i,u_i,s_i\}$. This subalgebra is not closed under
the crossing symmetry and so the crossing symmetry implies that
the $R$-matrix cannot be written as a polynomial in the braid matrix.

This $R$-matrix can also be written in terms of the basis of
orthogonal idempotents in (\ref{oidem}) as follows:
\begin{multline*}
-u^2q^2(q^2-1)^2R_i(u) = \\
(u-q^3)(u+q)(uq-1)(uq^3+1)\left[\frac1{2\delta} (u_i-t_i)\right]\\
(u-q^3)(u+q)(uq-1)(q^3+u)\left[\frac12 (f_i+s_i)\right] \\
(u-q^3)(u+q)(q-u)(q^3+u)\left[ e_i-\frac1\delta u_i\right]\\
(u-q^3)(1+uq)(q-u)(q^3+u)\left[\frac12 (f_i-s_i)\right] \\
(1-uq^3)(1+uq)(q-u)(q^3+u)\left[ \frac1{2\delta}(u_i+t_i)\right] 
\end{multline*}
In particular, the idempotents are independent of the spectral parameter
and so this $R$-matrix also has the property that
\begin{equation}\label{com}
R_i(u)R_i(v) = R_i(v)R_i(u)
\end{equation}
for all $u$ and $v$.

\bibliographystyle{halpha}
\bibliography{dilute}

\end{document}